\documentclass{amsart}
\usepackage{amsmath, amsfonts, amssymb}
\newcommand{\dal}{\square}
\newcommand{\R}{{\mathbb R}}
\newcommand{\ve}{\varepsilon}
\newcommand{\pa}{\partial}

\newcommand{\jb}[1]{\left\langle #1 \right\rangle}
\newcommand{\Dpc}{D_{+,c}}
\newcommand{\Dmc}{D_{-,c}}
\newcommand{\Dpmc}{D_{\pm, c}}

\newtheorem{theorem}{Theorem}
\newtheorem{lemma}{Lemma}

\numberwithin{equation}{section}
\numberwithin{theorem}{section}
\numberwithin{lemma}{section}
\title[Decay estimates of a tangential derivative]{%
Decay estimates of a tangential derivative to 
the light cone for the wave equation and their application}
\author[S.~Katayama]{Soichiro Katayama}
\address{Department of Mathematics, Wakayama University, 930 Sakaedani, Wakayama 640-8510, Japan}
\email{katayama@center.wakayama-u.ac.jp}
\author[H.~Kubo]{Hideo Kubo}
\address{Department of Mathematics, Graduate School of Science, Osaka University,
Toyonaka, Osaka 560-0043, Japan}
\email{kubo@math.sci.osaka-u.ac.jp}
\thanks{The first and the second author were
partially supported by Grant-in-Aid for Young Scientists (B)
(No.~16740094), MEXT, and by  Grant-in-Aid for Science Research (No.17540157), JSPS, respectively}
\thanks{Published in SIAM Journal on Mathematical Analysis Vol.~{\bf 39} (2008), no.~6, 1851--1862.}
\begin{document}
\begin{abstract}
We consider wave equations in three space dimensions and obtain
new weighted {$L^\infty$-$L^\infty$} estimates for a
tangential derivative to the light cone. As an application, we
give a new proof of the global existence theorem, which was
originally proved by Klainerman and Christodoulou, for systems of
nonlinear wave equations under the null condition. Our new proof
has the advantage of using neither the scaling nor the Lorentz
boost operators.   
\end{abstract}
\maketitle
\section{Introduction}\label{sec1}
Solutions to the Cauchy problem for nonlinear wave equations with
quadratic nonlinearity in three space dimensions may blow up in
finite time no matter how small initial data are, and we have to
impose some special condition on the nonlinearity to get global
solutions. The null condition is one of such conditions and is
associated with the null forms $Q_0$ and $Q_{ab}$, which are
given~by
\begin{align}
Q_0(v,w;c)  = & (\pa_t v)(\pa_t w)-c^2 (\nabla_x v)\cdot (\nabla_x w),\\
Q_{ab}(v,w)  = & (\pa_a v)(\pa_b w)-(\pa_b v)(\pa_a w)
\quad \mbox{($0\le a<b\le 3$)}
\end{align}
for $v=v(t,x)$ and $w=w(t,x)$, where $c$ is a positive constant
corresponding to the propagation speed, $\pa_0=\pa_t=\pa/\pa t$,
and $\pa_j=\pa/\pa x_j$ ($j=1, 2, 3$). More precisely, let $c>0$
and consider the Cauchy problem for
\begin{equation}
\label{WaveEq} \dal_c u_i=F_i(u, \pa u, \nabla_x\pa u)\quad
\mbox{in $(0,\infty)\times \R^3$\quad {\rm ($1\le i\le m$)}}
\end{equation}
with initial data
\begin{equation}
\mbox{$u=\ve f$ and $\pa_t u=\ve g$\quad at $t=0$,}
\label{InitialData}
\end{equation}
where $\dal_c=\pa_t^2-c^2\Delta_x$, $u=(u_j)$, $\pa u=(\pa_a
u_j)$, and $\nabla_x\pa u=(\pa_k\pa_a u_j)$ with $1\le j\le m$,
$1\le k\le 3$, and $0\le a\le 3$, while $\ve$ is a positive
parameter. Let $F=(F_i)_{1\le i\le m}$ be quadratic around the
origin in its arguments and the system be quasi-linear. In other
words, we assume that each $F_i$ has the form
\begin{equation}
F_i(u, \pa u, \nabla_x\pa u)= \sum_{{1\le j\le m \atop 1\le k\le
3,\ 0\le a\le 3}}c_{ka}^{ij}(u, \pa u)\pa_k\pa_a u_j+d_i(u, \pa
u), \label{quasi}
\end{equation}
where $c_{ka}^{ij}(u,\pa u)=O(|u|+|\pa u|)$ and $d_i(u, \pa
u)=O(|u|^2+|\pa u|^2)$ around $(u, \pa u)=(0,0)$. Without loss of
generality, we may assume $c^{ij}_{k\ell}=c^{ij}_{\ell k}$ for
$1\le i, j\le m$ and $1\le k, \ell \le 3$. In addition, we always
assume the symmetry condition
$$
c^{ij}_{ka}=c^{ji}_{k a}\quad \mbox{for $1\le i,j\le m$,\quad
$1\le k \le 3$, and $0\le a\le 3$}.
$$
Then it is well known that the null condition (for the above
system (\ref{WaveEq})) is satisfied if and only if the quadratic
terms of $F_i$ ($1\le i\le m$) can be written as linear
combinations of the null forms $Q_0(u_j, \pa^\alpha u_k;c)$ and
$Q_{ab}(u_j, \pa^\alpha u_k)$ with $1\le j, k\le m$, $0\le a<b\le
3$, and $|\alpha|\le 1$, where
$\pa^\alpha=\pa_0^{\alpha_0}\pa_1^{\alpha_1}\pa_2^{\alpha_2}\pa_3^{\alpha_3}$
for a multi-index $\alpha=(\alpha_0, \alpha_1, \alpha_2,
\alpha_3)$ (refer to~\cite{Chr86} and~\cite{Kla86} for the precise
description of the null condition). Klainerman~\cite{Kla86} and
Christodoulou~\cite{Chr86} proved the following global existence
theorem independently by different methods.

\begin{theorem}[Klainerman \cite{Kla86}, Christodoulou \cite{Chr86}]
\label{KlaThm} Suppose that the null condition is satisfied. Then,
for any $f$, $g\in C_0^\infty(\R^3;\R^m)$, there exists a positive
constant $\ve_0$ such that the Cauchy problem {\rm
(\ref{WaveEq})--(\ref{InitialData})} admits a unique global
solution $u\in C^\infty([0, \infty)\times \R^3; \R^m)$ for any
$\ve\in (0, \ve_0]$.
\end{theorem}

Christodoulou used the so-called conformal method which is based
on Penrose's conformal compactification of Minkowski space. On the
other hand, Klainerman used the vector field method and showed the
above theorem by deriving some decay estimates in the original
coordinates. In Klainerman's proof, he introduced vector fields
$$
 L_{c,j}=\frac{x_j}{c}\pa_t+ct\pa_j\quad (1\le j\le 3),\quad
 \Omega_{ij}=x_i\pa_j-x_j \pa_i\quad (1\le i<j\le 3),
$$
which are the generators of the Lorentz group, and the scaling
operator
$$
S=t\pa_t+x\cdot \nabla_x.
$$
These vector fields play an important role in getting Klainerman's
weighted $L^1$-$L^\infty$ estimates for wave equations (see also
H{\"o}rmander~\cite{Hoe88}). In addition, using them, we can see
that an extra decay factor is expected from the null forms. For
example, we have
\begin{align}
 \label{KlaEx}
Q_0(v, w; c) =& \frac{1}{t+r} \biggl\{ (\pa_t v)\bigl(S w+c L_{c,
r} w \bigr)
{}-c\sum_{j=1}^3 (L_{c,j} v)(\pa_j w) \\
 & \qquad\qquad {}-c^2(Sv)(\pa_r w) {}+c^2 \sum_{j\ne k}
\omega_k(\Omega_{jk} v)(\pa_j w) \biggr\}, \nonumber
\end{align}
where $r=|x|$, $\omega=(\omega_1,\omega_2, \omega_3)=x/r$,
$\pa_r=\sum_{j=1}^3 \omega_j \pa_j$, $L_{c,r} =\sum_{j=1}^3
\omega_j L_{c,j}$, and $\Omega_{ij}=-\Omega_{ji}$ for $1\le j<i\le
3$.

Among the above vector fields, the Lorentz boost fields $L_{c,j}$
depend on the propagation speed $c$, and they are unfavorable when
we consider the multiple speed case. Thus, the vector field method
without the Lorentz boost fields was developed by many authors
(see Kovalyov~\cite{Kov87,Kov89}, Klainerman and
Sideris~\cite{Kla-Sid96}, Yokoyama~\cite{Yok00}, Kubota and
Yokoyama~\cite{Kub-Yok01}, Sideris and Tu~\cite{Sid-Tu01},
Sogge~\cite{Sog03}, Hidano~\cite{Hid04},
Katayama~\cite{Kat04:02,Kat07p}, and Katayama and
Yokoyama~\cite{Kat-Yok06}, for example). In place of
(\ref{KlaEx}), the following identity was used in the above works
relating to the null condition for the multiple speed case:
\begin{align}
\label{YokoEx}
Q_0(v,w; c)
 =& \frac{1}{t^2}
(Sv+(ct-r)\pa_r v)(Sw-(ct+r)\pa_r w) \\*
 &{}+\frac{c}{t}
\left\{(Sv)(\pa_r w)-(\pa_r v)(Sw)\right\}
{}+\frac{c^2}{r}\sum_{j\ne k} {\omega_k}(\pa_jv)(\Omega_{jk}w),
\nonumber
\end{align}
whose variant was introduced by Hoshiga and Kubo~\cite{Hos-Kub00}.
Equation~(\ref{YokoEx}) leads to a good estimate in the region
$r>\delta t$ with some small $\delta>0$, because $r$ is equivalent
to $t+r$ in this region. Note that the operator $S$ is still used
in (\ref{YokoEx}), and this is the only reason why $S$ was adopted
in~\cite{Kat04:02,Kub-Yok01,Yok00}, because these works are based
on variants of $L^\infty$-$L^\infty$ estimates due to
John~\cite{Joh79} and Kovalyov~\cite{Kov87}, where only $\pa_a$
and $\Omega_{ij}$ are used (see Lemma~\ref{InhomDecay} below).

Our aim here is to get rid of not only $L_{c,j}$, but also $S$
from the estimate of the null forms, and prove
Theorem~\ref{KlaThm} using only $\pa_a$ and $\Omega_{jk}$. Though
the usage of the scaling operator $S$ has not caused any serious
difficulty in the study of the Cauchy problem for nonlinear wave
equations so far, we believe that it is worthwhile developing a
simple approach with a smaller set of vector fields. For this
purpose, we make use of the identity
\begin{align}
\label{KataEx}Q_0(v,w; c)
=&\frac{1}{2}\bigl\{(\Dpc v)(\Dmc w)+(\Dmc v)(\Dpc w)\bigr\} \\
 & {}+\frac{c^2}{r}\sum_{j\ne k} {\omega_k}(\pa_jv)(\Omega_{jk}w),
\nonumber
\end{align}
where $\Dpmc=\pa_t\pm c\pa_r$. Note that this identity was already
used implicitly to obtain identities like (\ref{YokoEx})
(see~\cite{Sid-Tu01}, for example). In view of (\ref{KataEx}),
what we need to treat the null forms is an enhanced decay estimate
for the tangential derivative $D_{+,c}$ to the light cone. We can
say that, in the previous works, this enhanced decay has been
observed through
$$
D_{+,c}=\frac{1}{t}\bigl(S+(ct-r)\pa_r\bigr)\ \mbox{or }
D_{+,c}=\frac{1}{ct+r}\bigl(cS+c L_{c,r}\bigr)
$$
with the help of $S$ or also $L_{c,r}=\sum_{j=1}^3\omega_j L_{c,j}$.

In this paper, we take a different approach. We will establish the
enhanced decay of $\Dpc u$ for the solution $u$ to the wave
equation directly. We formulate it as a weighted
$L^\infty$-$L^\infty$ estimate in Theorem~\ref{D+Es} below, which
is our main ingredient in this paper. The point is that such an
estimate can be derived by using only $\pa_a$ and $\Omega_{ij}$.
This type of approach to $\Dpc$ goes back to the work of
John~\cite{Joh83}.
\section{The Main Result}\label{sec2}

Before stating our result precisely, we introduce several
notations. We put $Z=\{Z_a\}_{1\le a\le 7}=\{(\pa_a)_{0\le a \le
3}, (\Omega_{jk})_{1\le j<k\le 3}\}$. For a multi-index
$\alpha=(\alpha_1, \ldots, \alpha_7)$, we define $Z^\alpha=
Z_1^{\alpha_1}Z_2^{\alpha_2}\cdots Z_7^{\alpha_7}$. For a function
$v=v(t,x)$ and a nonnegative integer $s$, we define
\begin{equation}
 |v(t,x)|_s=\sum_{|\alpha|\le s} |Z^\alpha v(t,x)|\ \mbox{and }
 \|v(t,\cdot)\|_{s}=\bigl\||v(t,\cdot)|_s \bigr\|_{L^2(\R^3)}.
\end{equation}

We put $\jb{a}=\sqrt{1+a^2}$ for $a\in \R$. Let $c$ be a positive
constant, and we fix arbitrary positive constants $c_j$ ($1\le
j\le N$) (our theorem is true for any choice of these constants
$c_j$, but when we apply our estimate to nonlinear problems, we
usually choose $c_j$ as the propagation speeds and $N$ as the
number of different propagation speeds in the system; $c$ is also
chosen from these propagation speeds). We define
\begin{equation}
 w(t,r)=w(t,r; c_1,\ldots, c_N)=\min_{0\le j\le N} \jb{c_jt-r}
 \label{Weight}
\end{equation}
with $c_0=0$, and we define
\begin{equation}
A_{\rho, \mu, s}[G;c](t,x)=\sup_{(\tau, y)\in \Lambda_{c}(t,x)}
 |y| \jb{\tau+|y|}^\rho w(\tau, |y|)^{1+\mu} |G(\tau, y)|_s
\label{InhomWei}
\end{equation}
for $\rho, \mu \ge 0$, a nonnegative integer $s$, and a smooth
function $G=G(t,x)$, where $\Lambda_{c}(t,x) =\{(\tau,y)\in
[0,t]\times \R^3\,;\, |y-x|\le c(t-\tau) \}$. We also define
\begin{equation}
B_{\rho, s}[\phi, \psi;c](t,x)=\sup_{y\in \Lambda'_{c}(t,x)} \jb{|y|}^{\rho}
\bigl(|\phi(y)|_{s+1}+|\psi(y)|_s\bigr)
\label{HomWei}
\end{equation}
for $\rho\ge 0$, a nonnegative integer $s$, and smooth functions
$\phi$ and $\psi$ on $\R^3$, where $\Lambda'_{c}(t,x)=\{y\in
\R^3\,;\, |y-x|\le ct \}$.

The following theorem is our main result.

\begin{theorem}\label{D+Es}
Assume $1\le \kappa\le 2$ and $\mu>0$.

\begin{enumerate}
\item[{\rm (i)}] Let $u$ be the solution to
$$
\dal_c u=G\quad \mbox{in $(0,\infty)\times \R^3$}
$$
with initial data $u=\pa_t u=0$ at $t=0$. Then there exists a
positive constant $C$, depending on $\kappa$ and $\mu$, such that
\begin{eqnarray}
\label{D+Es-01}
& & \jb{|x|} \jb{t+|x|} \jb{ct-|x|}^{\kappa-1}
\left\{\log(2+t+|x|)\right\}^{-1}
|\Dpc u(t,x)| \\
& & \qquad\qquad\qquad\qquad\qquad\qquad\qquad\qquad\qquad
\le C A_{\kappa, \mu, 2}[G; c](t,x)
\nonumber
\end{eqnarray}
for $(t,x)\in (0,\infty)\times \R^3$ with $x\ne 0$, where
$A_{\kappa, \mu, 2}$ is given by $(\ref{InhomWei})$.

\hspace*{1.5pc}Moreover, if $1<\kappa< 2$, then for any
$\delta>0$, there exists a constant C, depending on $\kappa$,
$\mu$, and $\delta$, such that
\begin{equation}
\label{D+Es-01'}
\jb{t+|x|}^2 \jb{ct-|x|}^{\kappa-1}
|\Dpc u(t,x)|
\le C A_{\kappa, \mu, 2}[G; c](t,x)
\end{equation}
for $(t,x)\in (0,\infty)\times \R^3$ satisfying $|x|>\delta t$.

\item[{\rm (ii)}] Let $u^*$ be the solution to
$$
\dal_c u^*=0 \mbox{ in $(0,\infty)\times \R^3$}
$$
with initial data $u^*=\phi$ and $\pa_t u^*=\psi$ at $t=0$.
Then we have
\begin{equation}
\label{D+Es-02}
\jb{|x|}\jb{t+|x|} \jb{ct-|x|}^{\kappa-1}|\Dpc u^*(t,x)|
\le C B_{\kappa+\mu+1,2}[\phi, \psi; c](t,x)
\end{equation}
for $(t,x)\in (0,\infty)\times \R^3$ with $x\ne 0$,
where $B_{\kappa+\mu+1, 2}$ is given by $(\ref{HomWei})$.
\end{enumerate}
\end{theorem}

{\it Remark}.
(1) Similar estimates for radially symmetric solutions
are obtained by Katayama \cite{Kat07p}.
\smallskip\\
(2) Suppose that $A_{\kappa,\mu, 2}[G; c](t,x)$ is bounded on
$[0,\infty)\times\R^3$ for some $\kappa \in[1,2)$ and $\mu>0$ and
that $u$ solves $\dal_cu=G$ with zero initial data. Then, from
Lemma~\ref{InhomDecay} below, we see that $u$ and $\pa u$ decay
like $\jb{t}^{-1}\Psi_{\kappa-1}(t)$ along the light cone
$ct=|x|$, where $\Psi_\rho(t)=\log (2+t)$ if $\rho=0$, and
$\Psi_\rho(t)=1$ if $\rho>0$. Compared with this decay rate, we
find from (\ref{D+Es-01}) and (\ref{D+Es-01'}) that $\Dpc u$ gains
extra decay of $\jb{t}^{-1}$ and behaves like
$\jb{t}^{-2}\Psi_{\kappa-1}(t)$ along the light cone.
\smallskip\\
(3) For tangential derivatives $T_{c,j}=(x_j/|x|)\pa_t+c\pa_j$ ($1\le j\le 3$),
Alinhac showed~that
$$
\left(\int_0^t \int_{\R^3} \bigl(1+\bigl|c\tau-|x|\,\bigr|\bigr)^{-\rho} |T_{c,j} u(\tau, x)|^2 dx d\tau \right)^{1/2}
$$
with $\rho>1$ is bounded by $\|\pa u(0,
\cdot)\|_{L^2(\R^3)}+\int_0^t\|\dal_c
u(\tau,\cdot)\|_{L^2(\R^3)}d\tau$ (see~\cite{Ali04}, for example).
Observe that $T_{c,j}$ is closely connected to $\Dpc$. In fact, we
have $\Dpc=\sum_{j=1}^3(x_j/|x|)T_{c,j}$. Though Alinhac's
estimate does not need $S$ and means enhanced decay of tangential
derivatives implicitly, it seems difficult to recover a pointwise
decay estimate from his weighted space-time estimate. On the other
hand, Sideris and Thomases \cite{Sid-Tho06} obtained the estimate
for $\left\|\bigl(1+\bigl|ct+|\cdot|\,\bigr|\bigr)T_{c,
j}u(t,\cdot) \right\|_{L^2(\R^3)}$; however, $S$ is used in their
estimate.
\smallskip\\
(4) The exterior problem for systems of nonlinear wave equations with the single or multiple
speed(s) is also widely studied (see Metcalfe, Nakamura, and
Sogge~\cite{MetNaSo05b} and Metcalfe and Sogge~\cite{MetSo07} and
the references cited therein). In the exterior domains, because of
their unbounded coefficients on the boundary, the Lorentz boosts
are unlikely to be applicable even for the single speed case. This
is another reason why the vector field method without the Lorentz
boosts is widely studied. In addition, $S$ also causes a technical
difficulty in the exterior problems. We will discuss the exterior
problem in a subsequent paper, and we will not go into further
details here.
\medskip

We will prove Theorem~\ref{D+Es} in the next section, after
stating some known weighted $L^\infty$-$L^\infty$ estimates for
wave equations. Though we can apply our theorem to exclude $S$
from the proof of the multiple speed version of
Theorem~\ref{KlaThm} in~\cite{Kat04:02,Kub-Yok01,Yok00}, we
concentrate on the single speed case for simplicity, and we will
give a new proof, without using $S$ and $L_{c,j}$, of
Theorem~\ref{KlaThm} in section~\ref{sec4} as an application of
our main theorem.

Throughout this paper, various positive constants, which may
change line by line, are denoted just by the same letter~$C$.

\section{Proof of Theorem \ref{D+Es}}
\label{sec3}

For $c>0$, $\phi=\phi(x)$, and $\psi=\psi(x)$, we write
$U_c^*[\phi, \psi]$ for the solution $u$ to the homogeneous wave
equation $\dal_c u=0$ in $(0,\infty) \times \R^3$ with initial
data $u=\phi$ and $\pa_t u=\psi$ at $t=0$. Similarly, for $c>0$
and $G=G(t,x)$, we write $U_c[G]$ for the solution $u$ to the
inhomogeneous wave equation $\dal_c u=G$ in $(0, \infty)\times
\R^3$ with initial data $u=\pa_t u=0$ at $t=0$.

For $U_c^*[\phi, \psi]$ we have the following.

\begin{lemma}\label{HomDecay}
Let $c>0$. Then, for $\kappa> 1$, we have
\begin{align}
\label{AsakuraEs} & \jb{t+|x|}\jb{ct-|x|}^{\kappa-1} |U_c^*[\phi, \psi](t,x)| \\
 & \qquad\qquad \le C
\sup_{y\in \Lambda'_c(t,x)} \jb{|y|}^{\kappa}\left(
\jb{|y|}|\phi(y)|_{1}+|y|\,|\psi(y)|
\right)
\nonumber
\end{align}
for $(t,x)\in [0,\infty)\times \R^3$.
\end{lemma}

For the proof, see Katayama and Yokoyama~\cite[Lemma
3.1]{Kat-Yok06} (see also Asakura~\cite{Asa86} and Kubota and
Yokoyama~\cite{Kub-Yok01}).

After the pioneering work of John \cite{Joh79}, a wide variety of
weighted $L^\infty$-$L^\infty$ estimates for $U_c[G]$ and $\pa
U_c[G]$ have been obtained
(see~\cite{Asa86,Kat04:02,Kat05,Kat-Mat05,Kat-Yok06,Kov87,Kov89,Kub-Yok01,Yok00}).
Here we restrict our attention to what will be used directly in
our proofs of Theorems~\ref{KlaThm} and~\ref{D+Es}.

\begin{lemma}\label{InhomDecay}
Let $c>0$. Define
\begin{align}
 \Phi_{\rho}(t,r)  =&
 \left\{
 \begin{array}{ll}
 \log\bigl(
 2+\jb{t+r}\jb{t-r}^{-1}
 \bigr) &\mbox{ if $\rho=0$},\\
 \jb{t-r}^{-\rho} & \mbox{ if $\rho>0$},
 \end{array}
 \right.
\label{KubYokWei}\\
\Psi_\rho(t) =&
 \left\{
 \begin{array}{ll}
 \log(2+t) &\mbox{ if $\rho=0$},\\
 1 & \mbox{ if $\rho>0$}.
 \end{array}
 \right.
\end{align}
Assume $\kappa\ge 1$ and $\mu>0$. Then we have
\begin{align}
\label{KatYok03}
 & \jb{t+|x|}\Phi_{\kappa-1}(ct, |x|)^{-1}
|U_c[G](t,x)|\le C A_{\kappa, \mu, 0}[G;c](t,x),\\
 & \jb{|x|}\jb{ct-|x|}^{\kappa}\Psi_{\kappa-1}(t)^{-1}|\pa U_c[G](t,x)|
 \le C A_{\kappa, \mu, 1}[G;c](t,x)
\label{KatYok05}
\end{align}
for $(t,x)\in [0,\infty)\times \R^3$, where $A_{\kappa, \mu,
s}[G;c]$ is given by $(\ref{InhomWei})$.
\end{lemma}

\begin{proof}
For the proof of (\ref{KatYok03}), see Katayama and
Yokoyama~\cite[equation~(3.6) in Lemma 3.2, and section
8]{Kat-Yok06} for $\kappa>1$ and Katayama~\cite{Kat07p} for
$\kappa=1$.

Next we consider (\ref{KatYok05}) with $\kappa>1$. From Lemma 8.2
in~\cite{Kat-Yok06}, we find that (\ref{KatYok05}) with $\pa
U_c[G]$ replaced by $U_c[\pa G]$ is true. Now (\ref{KatYok05})
follows immediately from Lemma~\ref{HomDecay}, because we have
$\pa_a U_c[G]=U_c[\pa_a G]+\delta_{a0}U_c^*[0, G(0, \cdot)]$ for
$0\le a\le 3$ with the Kronecker delta $\delta_{ab}$, and
$\jb{|y|}^{\kappa+1}|y|\,|G(0,y)|\le C A_{\kappa, \mu, 1}[G;
c](t)$ (note that we have $w(0,r)=\jb{r}$). Equation
(\ref{KatYok05}) for the case $\kappa=1$ can be treated similarly
(see \cite{Kub-Yok01} and~\cite{Kat04:02}). \qquad\end{proof}

Note that we will use (\ref{KatYok05}) in the proof of
Theorem~\ref{KlaThm} but not in that of Theorem~\ref{D+Es}.

Now we are in a position to prove Theorem~\ref{D+Es}. Suppose that
all the assumptions in Theorem~\ref{D+Es} are fulfilled. Without
loss of generality, we may assume $c=1$.

For simplicity of exposition, we write $D_\pm$ for $D_{\pm,
1}=\pa_t\pm \pa_r$. Similarly, $U^*[\phi, \psi]$, $U[G]$,
$A_{\rho, \mu, s}(t,x)$, and $B_{\rho,s}(t,x)$ denote $U_1^*[\phi,
\psi]$, $U_1[G]$, $A_{\rho, \mu, s}[G;1](t,x)$, and \break
$B_{\rho, s}[\phi, \psi;1](t,x)$, respectively.

First we prove (\ref{D+Es-01}). Assume $0<r=|x|\le 1$. We have
$$
 |D_+ u|\le |\pa_t u|+|\nabla_x u| \le \sum_{0\le a\le 3}
 \left|U[\pa_a G]\right|+\left|U^*[0, G(0, \cdot)]\right|.
$$
From (\ref{KatYok03}) in Lemma~\ref{InhomDecay}, we get
\begin{equation}
\jb{t+r}\Phi_{\kappa-1}(t,r)^{-1}
\left|U[\pa_a G](t,x)\right|
\le C A_{\kappa, \mu, 1}(t,x),
\label{B-01}
\end{equation}
while Lemma~\ref{HomDecay} leads to
\begin{align*}
\jb{t+r}\jb{t-r}^\kappa \left|U^*[0, G(0, \cdot)](t,x)\right|
 \le & C \sup_{y\in \Lambda'_1(t,x)} |y|\jb{|y|}^{\kappa+1}|G(0,y)|\\
 \le & C A_{\kappa, \mu, 0}(t,x).
\end{align*}
Thus we obtain (\ref{D+Es-01}) for $0<|x|\le 1$.

We set $v(t, r, \omega)=ru(t, r\omega)$ for $r>0$ and $\omega\in
S^2$. Then we have
\begin{equation}
 D_-D_+v(t,r,\omega)=rG(t,r\omega)
 {}+\frac{1}{r}\sum_{1\le j<k\le 3} \Omega_{jk}^2 u(t, r\omega).
 \label{rad}
\end{equation}
Let $r=|x|\ge 1$ and $1\le \kappa\le 2$. From (\ref{KatYok03}), we
get
\begin{align}
\label{Est01}
\frac{1}{r}\sum_{1\le j<k\le 3} |\Omega_{jk}^2 u(t, r\omega)|
 \le & C \jb{r}^{-1}\jb{t+r}^{-1}
\Phi_{\kappa-1}(t,r)A_{\kappa, \mu, 2}(t,r\omega) \\
 \le & C \jb{t+r}^{-\kappa}(\jb{r}^{-1}+\jb{t-r}^{-1})
A_{\kappa, \mu, 2}(t, r\omega),
\nonumber
\end{align}
where $\Phi_{\kappa-1}$ is from (\ref{KubYokWei}). It is easy to
see that
\begin{equation}
|rG(t, r\omega)|\le \jb{t+r}^{-\kappa} w(t, r)^{-1-\mu}
A_{\kappa, \mu, 0}(t,r\omega).
\label{Est02}
\end{equation}
Note that we have
$$
A_{\kappa, \mu,s}(\tau, (t+r-\tau)\omega)\le A_{\kappa, \mu, s}(t,
r\omega)\quad \mbox{for $0\le \tau\le t$}.
$$
Therefore, by (\ref{rad}), (\ref{Est01}), and (\ref{Est02}), we get
\begin{align}
\label{Est03}
|D_+v(t,r,\omega)|
 =& \left|\int_0^t \frac{d}{d\tau} (D_+v)(\tau, t+r-\tau, \omega)
d\tau\right| \\
 =& \left|\int_0^t (D_-D_+v)(\tau, t+r-\tau, \omega) d\tau\right|
\nonumber\\
 \le & C\jb{t+r}^{-\kappa} A_{\kappa, \mu, 2}(t,r\omega)\int_0^t \jb{t+r-\tau}^{-1}d\tau \nonumber\\
  & {}+C\jb{t+r}^{-\kappa} A_{\kappa, \mu, 2}(t,r\omega)\int_0^t \jb{t+r-2\tau}^{-1} d\tau
 \nonumber\\
  & {}+C\jb{t+r}^{-\kappa} A_{\kappa, \mu, 0}(t,r\omega)
 \int_0^t w(\tau, t+r-\tau)^{-1-\mu} d\tau
 \nonumber\\
 \le & C\jb{t+r}^{-\kappa} A_{\kappa, \mu, 2}(t, r\omega) \log(2+t+r).
 \nonumber
\end{align}
Since we have
$$
r D_+u(t,r\omega)= D_+v(t,r,\omega)-u(t, r\omega),
$$
from (\ref{Est03}) and (\ref{KatYok03}), we obtain
$$
\jb{r}\jb{t+r}\jb{t-r}^{\kappa-1}|D_+u(t,x)|
\le C \log(2+t+|x|) A_{\kappa, \mu, 2}(t,x)
$$
for $r=|x|\ge 1$. This completes the proof of (\ref{D+Es-01}).

To prove (\ref{D+Es-01'}), we first note that $\jb{t+r}\le
C\jb{r}$ for $r>\delta t$. Let $1<\kappa<2$. By the first line of
(\ref{Est01}), we have
\begin{equation}
\label{Est01'}
\frac{1}{r}\sum_{1\le j<k\le 3} |\Omega_{jk}^2 u(t, r\omega)|
\le C \jb{t+r}^{-2}\jb{t-r}^{-\kappa+1}A_{\kappa, \mu, 2}(t,r\omega)
\end{equation}
for $r>\max\{\delta t,1\}$. Obviously $r>\max\{\delta t,1\}$
yields $ t+r-\tau>\max\{\delta \tau,1\}$ for $0\le \tau \le t$.
Hence following similar lines to (\ref{Est03}), we obtain
$$
 |D_+v(t,r,\omega)|\le C \jb{t+r}^{-\kappa}
 A_{\kappa,\mu, 2}(t, r\omega) \quad\mbox{ for $r\ge \max\{\delta t, 1\}$}.
$$
This immediately implies (\ref{D+Es-01'}), because we already know
that $|D_+u|$ (resp.,\ $|D_+u-r^{-1}D_+v|$) has the desired bound
for $(\delta t<)r\le 1$ (resp.,\ $r\ge \max\{\delta t, 1\}$).

Now we are going to prove (\ref{D+Es-02}). Lemma~\ref{HomDecay}
immediately implies
$$
\jb{t+|x|}\jb{t-|x|}^{\kappa+\mu-1}|D_+ u^*(t,x)|
\le C B_{\kappa+\mu+1, 1}(t,x),
$$
which is better than (\ref{D+Es-02}) for $0<|x|\le 1$.
Lemma~\ref{HomDecay} also implies
\begin{align}
\label{C-02}
 & \frac{1}{r}\sum_{1\le j<k\le 3}|\Omega_{jk}^2 u^*(t,x)| \\
 & \qquad \le C \jb{r}^{-1}
\jb{t+r}^{-1}\jb{t-r}^{1-\kappa-\mu} B_{\kappa+\mu+1,2}(t,x)
\nonumber 
\\
 & \qquad \le C
\jb{t+r}^{-\kappa}(\jb{r}^{-1-\mu}+\jb{t-r}^{-1-\mu})B_{\kappa+\mu+1,2}(t,x)
\nonumber
\end{align}
for $r=|x|\ge 1$. Set $v^*(t, r, \omega)=ru^*(t,r\omega)$ for
$r\ge 0$ and $\omega \in S^2$. For $r\ge 1$, similarly to
(\ref{Est03}), we get
\begin{align*}
|D_+v^*(t,r,\omega)| = & \left|(D_+v^*)(0, t+r, \omega)
{}+\int_0^t (D_-D_+v^*)(\tau, t+r-\tau, \omega) d\tau\right| \\
 \le &
C\jb{t+r}^{-\kappa} B_{\kappa+1, 0}(t,r\omega)
\\
 & {}+C\jb{t+r}^{-\kappa} B_{\kappa+\mu+1, 2}(t, r\omega)
\int_0^t \jb{t+r-\tau}^{-1-\mu} d\tau
\\
 & {}+C\jb{t+r}^{-\kappa} B_{\kappa+\mu+1, 2}(t, r\omega)\int_0^t \jb{t+r-2\tau}^{-1-\mu} d\tau
\\
 \le & C\jb{t+r}^{-\kappa} B_{\kappa+\mu+1, 2}(t, r\omega),
\end{align*}
which ends up with
$$
\jb{r}\jb{t+r}\jb{t-r}^{\kappa-1}|D_+u^*(t,x)| \le C
B_{\kappa+\mu+1, 2}(t,x)
$$
for $r=|x|\ge 1$. This completes the proof of (\ref{D+Es-02}).
\qquad\qed
\section{Proof of Theorem \ref{KlaThm}}
\label{sec4}

As an application of Theorem~\ref{D+Es}, we give a new proof of
Theorem~\ref{KlaThm}. First we derive estimates for the null
forms.

\begin{lemma}\label{NFEs}
Let $c$ be a positive constant, and $v=(v_1, \ldots, v_M)$.
Suppose that $Q$ is one of the null forms. Then, for a nonnegative
integer $s$, there exists a positive constant $C_s$, depending
only on $c$ and $s$, such that
\begin{align*}
|Q(v_j, v_k)|_s \le & C_s \biggl\{ |\pa v|_{[s/2]}
\sum_{|\alpha|\le s} |\Dpc Z^\alpha v|
{}+|\pa v|_{s} \sum_{|\alpha|\le [s/2]} |\Dpc Z^\alpha v|\\
 & \qquad\quad\qquad\qquad\qquad\qquad
 {}+\frac{1}{r}\bigl(|\pa v|_{[s/2]}|v|_{s+1}+|v|_{[s/2]+1}|\pa v|_s
 \bigr)
\biggr\}.
\end{align*}
\end{lemma}

\begin{proof}
The case $Q=Q_0$ and $s=0$ follows immediately from
(\ref{KataEx}). We can obtain similar identities for other null
forms by using
$$
 (\pa_t, \nabla_x)=
\left(\frac{1}{2}, -\frac{x}{2cr}\right)\Dmc+
\left(\frac{1}{2}, \frac{x}{2cr}\right)\Dpc
{}-\left(0, \frac{x}{r^2} \wedge \Omega \right)
$$
with $\Omega=(\Omega_{23}, -\Omega_{13}, \Omega_{12})$ (see (5.2)
in Sideris and Tu \cite[Lemma 5.1]{Sid-Tu01}), and we can show the
desired estimate for $s=0$. Since $Z^\alpha Q(v_j, v_k)$ can be
written in terms of $Q_0(Z^\beta v_j, Z^\gamma v_k; c)$ and
$Q_{ab}(Z^\beta v_j, Z^\gamma v_k)$ ($0\le a<b\le 3$) with
$|\beta|+|\gamma|\le |\alpha|$, the desired estimate for general
$s$ follows immediately. \qquad\end{proof}

Now we are going to prove Theorem~\ref{KlaThm}. Without loss of
generality, we may assume $c=1$. Assume that the assumptions in
Theorem~\ref{KlaThm} are fulfilled. Let $u$ be the solution to
(\ref{WaveEq})--(\ref{InitialData}) on $[0,T)\times \R^3$, and we
set
\begin{align*}
e_{\rho, k}(t,x) =& \jb{t+|x|}\jb{t-|x|}^\rho
|u(t,x)|_{k+2}
{}+\jb{|x|}\jb{t-|x|}^{\rho+1}|\pa u(t,x)|_{k+1}\\
 & {}+\chi(t,x)\jb{t+|x|}^2\jb{t-|x|}^\rho\sum_{|\alpha|\le k}
|D_{+,1}Z^\alpha u(t,x)|
\end{align*}
for $\rho>0$ and a positive integer $k$, where $\chi(t,x)=1$ if
$|x|>(1+t)/2$, while $\chi(t,x)=0$ if $|x|\le (1+t)/2$. We fix
$\rho\in(1/2, 1)$ and $s\ge 8$, and assume that
\begin{equation}
\sup_{0\le t<T} \|e_{\rho, s}(t,\cdot)\|_{L^\infty(\R^3)} \le M\ve
\label{InductiveAs}
\end{equation}
holds for some large $M(>0)$ and small $\ve(>0)$, satisfying $M\ve
\le 1$. Our goal here is to get (\ref{InductiveAs}) with $M$
replaced by $M/2$. Once such an estimate is established, it is
well known that we can obtain Theorem~\ref{KlaThm} by the
so-called bootstrap (or continuity) argument.

In the following we always assume $M$ is large enough, and $\ve$
is sufficiently small. For simplicity of exposition, we will not
write dependence of nonlinearities on the unknowns explicitly.
Namely we abbreviate $F(u, \pa u, \nabla_x\pa u)(t,x)$ as
$F(t,x)$, and so~on.

First we evaluate the energy. For any nonnegative integer $k\le
2s$, (\ref{InductiveAs}) implies
\begin{equation}
 |F^{(2)}(t,x)|_{k}
 \le C M\ve \jb{|x|}^{-1}\jb{t-|x|}^{-1-\rho} |\pa u(t,x)|_{k+1},
\label{ene01}
\end{equation}
where $F^{(2)}$ denotes the quadratic terms of $F$. Put
$H=F-F^{(2)}$, and $Zu=(Z_1u, \ldots, Z_7 u)$. Since we have
\begin{equation}
\label{f01}
 \jb{r}^{-1}\jb{t-r}^{-1}\le C\jb{t+r}^{-1}\quad
\mbox{for any $(t,r)\in [0, \infty)\times [0, \infty)$,}
\end{equation}
and since $\jb{|x|}^{-1}|Z u|\le C|\pa u|$, from
(\ref{InductiveAs}) we obtain
\begin{align}
\label{ene02}
|H(t,x)|_{k}
 \le & C\bigl(|u|^3+|(u, \pa u)|_{[k/2]+1}^2
(|Zu|_{k-1}+|\pa u|_{k+1})\bigr) \\
 \le & CM^3\ve^3\jb{t+|x|}^{-3}\jb{t-|x|}^{-3\rho} \nonumber 
\\
 & {}+CM^2\ve^2 \jb{t+|x|}^{-1}\jb{t-|x|}^{-2\rho}|\pa u(t,x)|_{k+1}
\nonumber
\end{align}
for any nonnegative integer $k\le 2s$. Similarly to (\ref{ene01})
and (\ref{ene02}), using (\ref{f01}), we obtain
\begin{equation}
|F_{i,\alpha}(t,x)|\le CM\ve(1+t)^{-1}|\pa u(t,x)|_{2s}
{}+CM^3\ve^3
\jb{t+|x|}^{-3}\jb{t-|x|}^{-3\rho}
\label{ene02'}
\end{equation}
for $|\alpha|\le 2s$, where
$$
F_{i,\alpha}=Z^\alpha F_i-\sum_{j,k,a}c_{ka}^{ij}\pa_k\pa_a
(Z^\alpha u_j)
$$
with $c_{ka}^{ij}$ coming from (\ref{quasi}). It is easy to see
that
\begin{equation}
 \|\jb{t+|\cdot|}^{-3}\jb{t-|\cdot|}^{-3\rho}\|_{L^2(\R^3)}
 \le C (1+t)^{-2}
\end{equation}
for $\rho>1/2$. Therefore, from (\ref{ene02'}), we obtain
$$
\|F_{i,\alpha}(t,\cdot)\|_{L^2}
\le CM\ve (1+t)^{-1} \|\pa u(t,\cdot)\|_{2s}+CM^3\ve^3(1+t)^{-2}
$$
for $|\alpha|\le 2s$. We also have
$$
\sum_{j,k,a}|c_{ka}^{ij}(t,x)|_1\le CM\ve (1+t)^{-1}.
$$
Now, applying the energy inequality for the systems of perturbed
wave equations $\dal_1 (Z^\alpha u_i)-\sum_{j,k,a} c^{ij}_{ka}
\pa_k\pa_a (Z^\alpha u_j) =F_{i,\alpha}$, we find
$$
\frac{d}{dt}\|\pa u(t,\cdot)\|_{2s}\le
CM\ve(1+t)^{-1} \|\pa u(t, \cdot)\|_{2s}+CM^3\ve^3(1+t)^{-2},
$$
and the Gronwall lemma leads to
\begin{equation}
\|\pa u(t,\cdot)\|_{2s}\le C(\ve+M^3\ve^3) (1+t)^{C_0M\ve} \le
CM\ve (1+t)^{C_0M\ve} \label{Energy}
\end{equation}
with an appropriate positive constant $C_0$ which is independent
of $M$ (note that the energy inequality for the systems of
perturbed wave equations is available because of the symmetry
condition).

In the following, we repeatedly use Theorem~\ref{D+Es} and
Lemmas~\ref{HomDecay} and~\ref{InhomDecay} with the choice of
$N=1$ and $c_1=1(=c)$. In other words, from now on we put
$w(t,r)=\min\bigl\{\jb{r}, \jb{t-r}\bigr\}$. Note that we have
\begin{equation}
 \jb{r}^{-1}\jb{t-r}^{-1}\le C \jb{t+r}^{-1}w(t,r)^{-1},
\label{WeiIneq}
\end{equation}
which is more precise than (\ref{f01}).

By (\ref{Energy}) and the Sobolev-type inequality
$$
 \jb{|x|}|v(t,x)|\le C \|v(t,\cdot)\|_{2},
$$
whose proof can be found in Klainerman~\cite{Kla87}, we see that
\begin{equation}
\jb{|x|}|\pa u(t,x)|_{2s-2}\le CM\ve (1+t)^{C_0M\ve}.
\label{Highest}
\end{equation}
Using (\ref{WeiIneq}) and (\ref{Highest}), from (\ref{ene01}) and
(\ref{ene02}) with $k=2s-3$, we obtain
$$
|F(t,x)|_{2s-3} \le CM^2\ve^2 \jb{r}^{-1}\jb{t+|x|}^{-1}w(t,
|x|)^{-2\rho}(1+t)^{C_0M\ve},
$$
which implies
\begin{equation}
A_{1+\nu, 2\rho-1, 2s-3}[F;1](t,x)\le C M^2\ve^2
\jb{t+|x|}^{C_0M\ve+\nu},
\end{equation}
where $\nu$ is a positive constant to be fixed later (note that we
have $\jb{\tau+|y|}\le \jb{t+|x|}$ for $(\tau,y)\in
\Lambda_1(t,x)$). Since $2\rho>1$ and $1+\nu>1$, by
Lemmas~\ref{HomDecay} and~\ref{InhomDecay} with
Theorem~\ref{D+Es}, we obtain
\begin{eqnarray}
\label{Higher}
e_{0, 2s-5}(t,x)
& \le & e_{\nu, 2s-5}(t,x)\le C\ve+CM^2\ve^2\jb{t+|x|}^{C_0M\ve+\nu} \\
& \le &
CM\ve\jb{t+|x|}^{C_0M\ve+\nu}. \nonumber
\end{eqnarray}

Finally, we are going to estimate $e_{\rho, s}(t,x)$. By
(\ref{Higher}) and (\ref{ene01}) with $k=2s-6$, we have
$$
|F^{(2)}(t,x)|_{2s-6}
\le CM^2\ve^2 \jb{t+|x|}^{-2-\rho+C_0M\ve+\nu}\jb{|x|}^{-2}
$$
for $(t,x)$ satisfying $|x|\le (t+1)/2$. On the other hand,
(\ref{InductiveAs}), (\ref{Higher}), and Lemma~\ref{NFEs} imply
$$
|F^{(2)}(t,x)|_{2s-6} \le
CM^2\ve^2\jb{t+|x|}^{-3+C_0M\ve+\nu}\jb{t-|x|}^{-1-\rho}
$$
for $(t,x)$ satisfying $|x|\ge (t+1)/2$. Summing up, we obtain
\begin{equation}
\label{Lower01} |F^{(2)}(t,x)|_{2s-6} \le CM^2\ve^2
\jb{|x|}^{-1}\jb{t+|x|}^{-2+C_0M\ve+\nu} w(t,|x|)^{-1-\rho}.
\end{equation}
By the first line of (\ref{ene02}) with $k=2s-6$, using
(\ref{InductiveAs}) and (\ref{Higher}), we get
\begin{equation}
\label{Lower02}
|H(t,x)|_{2s-6}
\le CM^3\ve^3 \jb{|x|}^{-1}\jb{t+|x|}^{-2+C_0M\ve+\nu}
w(t,|x|)^{-2\rho}.
\end{equation}
Equations (\ref{Lower01}) and (\ref{Lower02}) yield
\begin{equation}
\label{Lower03}
|F(t,x)|_{2s-6}
\le
CM^2\ve^2 \jb{|x|}^{-1}\jb{t+|x|}^{-2+C_0M\ve+\nu} w(t,|x|)^{-2\rho}.
\end{equation}

Now we fix some $\nu$ satisfying $0<\nu<1-\rho$, and assume that
$\ve$ is sufficiently small to satisfy $-2+C_0M\ve+\nu\le
-1-\rho$. Then from (\ref{Lower03}) we find that
\begin{equation}
A_{1+\rho, 2\rho-1, 2s-6}[F;1](t,x)\le CM^2\ve^2.
\label{Lower04}
\end{equation}
Since we have $s+2\le 2s-6$, $1+\rho>1$, and $2\rho>1$, from
Theorem~\ref{D+Es}, Lemmas~\ref{HomDecay} and~\ref{InhomDecay}, we
obtain
\begin{equation}
e_{\rho, s}(t,x)\le C_1\bigl(\ve+M^2\ve^2\bigr)
\label{Lower05}
\end{equation}
for $(t,x)\in [0,T) \times \R^3$, with an appropriate positive
constant $C_1$ which is independent of $M$. Finally, if $M$ is
large enough to satisfy $4C_1\le M$, and $\ve$ is small enough to
satisfy $C_1M\ve \le 1/4$, by (\ref{Lower05}) we obtain
\begin{equation}
\sup_{0\le t<T} \|e_{\rho, s}(t,\cdot)\|_{L^\infty(\R^3)}\le \frac{M}{2}\ve,
\end{equation}
which is the desired result. This completes the proof.
\qquad\endproof


\end{document}